\documentclass[12pt,twoside]{article}

\usepackage[centertags]{amsmath}
\usepackage{amsfonts}
\usepackage{amssymb}
\usepackage{amsthm}
\usepackage{newlfont}

\newcommand{\RR}{\mathbb R}

\newcommand{\tr}{\operatorname{tr}}
\newcommand{\sign}{\operatorname{sign}}

\newcommand{\spa}{\operatorname{span}}

\newcommand{\image}{\operatorname{image}}

\newcommand{\End}{\operatorname{End}}

\newcommand{\grad}{\operatorname{grad}}
\newcommand{\rank}{\operatorname{rank}}
\newcommand{\diver}{\operatorname{div}}
\newcommand{\diag}{\operatorname{diag}}

\newtheorem{thm}{Theorem}[section]

\newtheorem{lem}[thm]{Lemma}
\newtheorem{prop}[thm]{Proposition}

\theoremstyle{definition}
\newtheorem{exam}[thm]{Example}

\theoremstyle{definition}
\newtheorem{defn}[thm]{Definition}
\theoremstyle{remark}
\newtheorem{rem}[thm]{Remark}
\title{Harmonic morphisms between degenerate semi-Riemannian manifolds}
\author{Alberto Pambira\thanks{University of Leeds, LS2 9JT, School of Mathematics; e-mail:
{\tt pambira@maths.leeds.ac.uk}}}
\date{}
\begin{document}
\maketitle
\begin{abstract}
In this paper we generalize harmonic maps and morphisms
to the \emph{degenerate semi-Riemannian category}, in the case
when the manifolds $M$ and $N$ are \emph{stationary} 
and the map $\phi :M\rightarrow N$ is \emph{radical-preserving}. We characterize geometrically the notion of
\emph{(generalized) horizontal (weak) conformality} and we obtain a
characterization for (generalized) harmonic morphisms in terms of (generalized) harmonic maps.
\end{abstract}
\section{Introduction and preliminaries}
\emph{Harmonic morphisms between (non-degenerate semi-)Riemannian manifolds} are
maps which preserve germs of harmonic functions. They are characterized in \cite{Fu78, Ish, {Fu96}}
as the subclass of harmonic maps which are \emph{horizontally weakly conformal}.
An up-to-date bibliography on this topic is given in \cite{G};
see also \cite{G-L-M} for a list of harmonic morphisms and construction techniques, 
and \cite{B-W} for a comprehensive account of the topic.

However, when the manifold $(M,g)$ is \emph{degenerate}, then it fails, in general,
to have a \emph{torsion-free, metric-compatible} connection; moreover,
in this case, the notion of `trace', with respect to the metric $g$,
does not make any sense, so that it is not possible to define the
`tension field' of a map, or, consequently, the notion of harmonic map, in the usual sense.

\emph{Degenerate manifolds} arise naturally in the semi-Riemannian category: for example the
restriction of a non-degenerate metric to a degenerate submanifold is a degenerate metric and
the Killing-Cartan form on a non-semi-simple Lie Group is a degenerate metric.

Such manifolds are playing an increasingly important role in quantum theory
and string theory, as the action and field equations of particles and strings often do not depend on the inverse metric and
are well-defined even when the metric becomes degenerate (cf. \cite{Ca-Ri}). For example,
an extension of Einstein's gravitational theory which contains degenerate metrics as possible solutions
might lead to space-times with no causal structure (cf. \cite{Ben}).

In the mathematical literature, degenerate manifolds have been studied under several names: \emph{singular 
Riemannian spaces}
(\cite{Moi,Vran,Vog1}), \emph{degenerate (pseudo- \emph{or} semi-Riemannian) manifolds} (\cite{Cr,Sul,K}),
\emph{lightlike manifolds} (\cite{D-B}), \emph{isotropic spaces} (\cite{Stru1,Stru2,Stru3,Stru4}),
\emph{isotropic manifolds} (\cite{Vog2}).

In this paper we define generalized harmonic maps and morphisms, characterize (generalized)
horizontally weakly conformal maps into four types (Theorem \ref{bigguy}), and give a
Fuglede--Ishihara-type characterization for generalized harmonic morphisms (Theorem \ref{bigguy2}).
We refer the reader to \cite{Pam} for further details.

In this section, we aim to introduce the necessary background on semi-Riemannian geometry which
will be used in the rest of the paper.
We shall assume that all vector spaces, manifolds etc. have finite dimension.
\subsection{Algebraic background}
Let $V$ be a vector space of dimension $m$.
\begin{defn}\label{metricsignature}
An \emph{inner product} on $V$ is a symmetric bilinear form
$\langle,\rangle =\langle ,\rangle _{V}$ on $V$. It is said to be \emph{non-degenerate} (on $V$) if $\langle w,w'\rangle =0$ for all
$w'\in V$ implies $w=0$, otherwise it is called \emph{degenerate}.

We shall refer to the pair $(V, \langle ,\rangle )$ as an \emph{inner product space}.
Given two subspaces
$W,W'\subseteq V$, we shall often write $W \perp _{V}W'$ to denote that $W$ is orthogonal to $W'$
(equivalently $W'$ is
orthogonal to $W$) with respect to the inner product $\langle ,\rangle _{V}$, i.e. $\langle w,w'\rangle =0$ for any
$w\in W$ and $w'\in W'$.
\end{defn}

Let $r,p,q \geq 0$ be integers and set $(\epsilon )_{ij}:=(\epsilon _{r,p,q})_{ij}$ to be the diagonal
matrix
\begin{displaymath}
(\epsilon )_{ij}=\diag (\underbrace{0, \ldots ,0}_{r - \textrm{times}},
\underbrace{-1, \ldots ,-1}_{p - \textrm{times}},\underbrace{+1, \ldots ,+1}_{q - \textrm{times}}).
\end{displaymath}
Given an inner product $\langle, \rangle$ on $V$, there exists a basis $\{ e_{i} \}$,
with $i=1, \ldots , m=\dim V$, of $V$ such that
$\langle e_{i},e_{i}\rangle =(\epsilon _{r,p,q})_{ij}$.
We call such a basis \emph{orthonormal}
and the triple $(r,p,q)$ is called the \emph{signature of the inner product
$\langle ,\rangle $}.
\begin{exam}
The \emph{standard $m$-Euclidean space $\RR ^{m}_{r,p,q}$ of signature $(r,p,q)$} is $\RR ^{m}$
endowed with the inner product $\langle ,\rangle _{r,p,q}$ defined by $\langle E_{i},E_{j}\rangle _{r,p,q}
:=(\epsilon _{r,p,q})_{ij}$; here $\{ E_{k} \} _{k=1}^{m}$ is the canonical basis $E_{1}=(1,0,\ldots ,0),
\ldots ,E_{m}=(0,\ldots ,0,1)$.
\end{exam}
\begin{defn}\label{frododeg}
A subspace $W$ of an inner product vector space $(V, \langle ,\rangle )$  is called \emph{degenerate}
(resp. \emph{null}) if there
exists a non-zero vector $X\in W$ such that $\langle X,Y \rangle =0$ for all $Y\in W$
(resp. if, for all $X,Y\in W$, we have
$\langle X,Y \rangle =0$). Otherwise $W$ is called \emph{non-degenerate} (resp. \emph{non-null}).
\end{defn}

Clearly if $W\neq \{ \mathbf{0} \}$ is null then it is degenerate. Moreover $W$ is degenerate
if and only if $\langle , \rangle |_{W}$ is degenerate, but this does not necessarily mean that $\langle ,\rangle $ 
is degenerate on $V$.

Given a vector space $V$, we define the \emph{radical of \/$V$} (cf. \cite{D-B}, p.1, 
\cite{K}, p.3 or \cite{ON}, p.53), denoted by ${\cal N}(V)$, to be the vector space:
\begin{displaymath}
{\cal N}(V):=V^{\perp}= \{ X\in V: \langle X,Y\rangle =0 \textrm{ for all } Y\in V \} .
\end{displaymath}

We notice (cf. \cite{ON}, p.49) that ${\cal N}(V)$ is a null subspace of $V$. Moreover, $V$ is non-degenerate if and only if
${\cal N}(V)= \{ \mathbf{0} \}$, and $V$ is null if and only if ${\cal N}(V)= V$.
Note that, for any subspace $W$ of $V$,
\begin{equation}\label{nig}
{\cal N}(V)\subseteq W^{\perp}.
\end{equation}

The following proposition generalizes two well-known facts of linear algebra (cf. \cite{ON}, chapter 2, Lemma 22).
\begin{prop}\label{dimensions}
For any subspace $W\subseteq V$ of an inner product space $(V,\langle ,\rangle )$ we have:
\begin{enumerate}
\item [(i)] $\dim W+\dim W^{\perp}=\dim V+ \dim ({\cal N}(V)\cap W)$;
\item [(ii)] $(W^{\perp})^{\perp}=W+{\cal N}(V)$.
\end{enumerate}
\end{prop}
\proof
Let $0 \leq t \leq \dim {\cal N}(V)$ be the integer $t=\dim {\cal N}(V)-\dim \big( W\cap {\cal N}(V)\big) $.
We can choose a basis $\{ e_{i} \} _{i=1}^{m}$ on $V$, `adapted' to ${\cal N}(V)$ and $W$, in the sense that
${\cal N}(V)=\spa (e_{1},\ldots ,e_{\dim {\cal N}(V)})$
and $W=\spa (e_{t+1}, \ldots ,e_{t+\dim W})$; claim \emph{(i)} follows immediately.

To prove \emph{(ii)} we note that
\begin{displaymath}
W+{\cal N}(V)\subseteq (W^{\perp})^{\perp}.
\end{displaymath}
From linear algebra (cf. \cite{T-T}, Theorem 1.9A) we have:
\begin{displaymath}
\dim (W+{\cal N}(V))=\dim W +\dim {\cal N}(V)-\dim (W\cap {\cal N}(V));
\end{displaymath}
on the other hand, by \emph{(i)} we get:
\begin{displaymath}
\dim W^{\perp}=\dim V+\dim (W\cap {\cal N}(V))- \dim W;
\end{displaymath}
on combining these and using (\ref{nig}) we obtain
\begin{displaymath}
\dim (W^{\perp})^{\perp} =\dim (W+{\cal N}(V));
\end{displaymath}
claim \emph{(ii)} follows.
\endproof

Let $W\subseteq V$ be a vector subspace of an inner product vector space $(V,\langle ,\rangle _{V})$ and
let $W^{\perp _{V}}$ be its orthogonal complement in $V$ with respect to $\langle ,\rangle _{V}$.
Denote by $\overline{V}$, $\overline{W}$ and $\overline{W^{\perp _{V}}}$ the spaces
\begin{equation}\label{biondabella}
\overline{V}:=V\big/ {\cal N}(V), \quad \overline{W}:=W\big/ ({\cal N}(V)\cap W), 
\textrm{ and } \, \overline{W^{\perp _{V}}}:=W^{\perp _{V}}\big/
{\cal N}(V),
\end{equation}
having noted that, by (\ref{nig}), ${\cal N}(V)\subseteq W^{\perp _{V}}$. Let us also denote by
$\langle ,\rangle _{\overline{V}}$ the inner product on $\overline{V}$ defined by
\begin{displaymath}
\langle \overline{v}, \overline{v'} \rangle _{\overline{V}}:=\langle v,v' \rangle _{V} \qquad (v,v'\in V).
\end{displaymath}
Note that this is well defined. For any subspace $E\subseteq \overline{V}$, let $E^{\perp _{\overline{V}}}$ denote its
orthogonal complement in $(\overline{V}, \langle ,\rangle _{\overline{V}})$. Then we have the following
\begin{prop}\label{bigdiscovery}
For any vector subspace $W\subseteq V$ we have the following canonical isomorphism:
\begin{equation}\label{sparaunnumero}
\overline{W}\cong (\overline{W^{\perp _{V}}})^{\perp _{\overline{V}}}.
\end{equation}
\end{prop}
\proof
Consider the composition
\begin{displaymath}
\theta :W\stackrel{i}{\hookrightarrow}V\stackrel{\pi _{V}}{\rightarrow}V\big/ {\cal N}(V)=:\overline{V},
\end{displaymath}
where $i:W\hookrightarrow V$ is the inclusion map and $\pi _{V}:V\rightarrow \overline{V}$ is the natural
projection. We have
\begin{displaymath}
\theta (W) \subseteq \left( W^{\perp _{V}}\big/ {\cal N}(V) \right) ^{\perp _{\overline{V}}};
\end{displaymath}
in fact, let $w\in W$ and $w'\in W^{\perp _{V}}$ and write $\theta (w):=\overline{w}$; then we have
\begin{displaymath}
0=\langle w,w'\rangle _{V}=\langle \overline{w},\overline{w}'\rangle _{\overline{V}} \, .
\end{displaymath}
Next, note that $\ker \theta ={\cal N}(V)\cap W$. In fact for any $w\in W$, we have
\begin{displaymath}
\theta (w)=0\iff \overline{w}=0 \iff w\in {\cal N}(V),
\end{displaymath}
and so the claim. Then $\theta$ factors to an injective map
\begin{displaymath}
\overline{\theta}:\overline{W}:=W\big/ {\cal N}(V)\cap W\longrightarrow
\left( W^{\perp _{V}}\big/ {\cal N}(V) \right) ^{\perp _{\overline{V}}}=:
(\overline{W^{\perp _{V}}})^{\perp _{\overline{V}}} \, .
\end{displaymath}

We show that this is an isomorphism, by calculating the dimension of the spaces on either side of the equation (\ref{sparaunnumero}).
On the left-hand side we have
\begin{displaymath}
\dim \overline{W}=\dim W- \dim ({\cal N}(V)\cap W);
\end{displaymath}
on the right-hand side, applying Proposition \ref{dimensions}, we get
\begin{displaymath}
\dim W^{\perp _{V}}= \dim V + \dim ({\cal N}(V)\cap W) -\dim W,
\end{displaymath}
so that
\begin{displaymath}
\dim \overline{W^{\perp _{V}}}= \dim V  + \dim ({\cal N}(V)\cap W)-\dim W -\dim {\cal N}(V)
\end{displaymath}
and, applying once more Proposition \ref{dimensions},
\begin{displaymath}
\begin{array}{rl}
\dim (\overline{W^{\perp _{V}}})^{\perp _{\overline{V}}}&=\dim \overline{V} 
-\big( \dim V  + \dim ({\cal N}(V)\cap W)-\dim W -\dim {\cal N}(V) \big) \\
                                                        &=\dim W-\dim ({\cal N}(V)\cap W) \\
                                                        &=\dim \overline{W},
\end{array}
\end{displaymath}
so that the map $\overline{\theta}$ is an isomorphism, and the claim follows.
\endproof
We shall use the Proposition above to identify
$\overline{W}$ and $(\overline{W^{\perp _{V}}})^{\perp _{\overline{V}}}$.
Thus, any subspace 
$K\subseteq \overline{W}$ will sometimes be considered as a subspace of $(\overline{W^{\perp _{V}}})^{\perp _{\overline{V}}}$
and \emph{vice versa}.
\subsection{Background on semi-Riemannian geometry}
\begin{defn}
Let $r,p,q$ be three non-negative integers such that $r+p+q=m$.
A \emph{semi-Riemannian metric $g$ of signature $(r,p,q)$ on an $m$-dimensional
smooth manifold $M$} is a smooth section of the symmetric square $\odot ^{2}T^{*}M$
which defines an inner product $\langle ,\rangle $ on each tangent space of constant
signature $(r,p,q)$. A \emph{semi-Riemannian manifold} is a pair $(M,g)$
where $M$ is a smooth manifold and $g$ is a semi-Riemannian metric
on $M$. When $r>0$ (resp. $r=0$, $r<m$, or $r=m$) $(M,g)$ is called \emph{degenerate}
(resp. \emph{non-degenerate}, \emph{non-null}, or \emph{null}).
\end{defn}

Let $\cal L$ denote the Lie derivative and let ${\cal N}={\cal N}(TM):=\cup _{x\in M}{\cal N}(T_{x}M)$;
${\cal N}$ is called the \emph{radical distribution on $M$}.
\begin{defn}[\cite{K}, Definition 3.1.3]\label{stationaryOK}
A semi-Riemannian manifold $(M,g)$ is said to be \emph{stationary} if
${\cal L}_{A}g=0$ for any locally defined
smooth section $A\in \Gamma ({\cal N})$.
\end{defn}

Such a manifold is also
called a \emph{Reinhart manifold} (cf. \cite{D-B}, p.49). The condition that $M$ be stationary is equivalent to
$\cal N$ being a Killing distribution (i.e. all vector fields in $\cal N$ are
Killing). Trivially a non-degenerate manifold is
stationary.

We introduce the following operator (\cite{K}, Definition 3.1.1):

\begin{defn}[Koszul derivative]\label{Koszul}
Let $(M,g)$ be a semi-Riemannian manifold. An
operator $D:\Gamma (TM) \times \Gamma (TM) \rightarrow \Gamma (TM)$ is
called a \emph{Koszul derivative on $(M,g)$} if, for any $X,Y,Z\in \Gamma
(TM)$, it satisfies the \emph{Koszul formula}
\begin{equation}\label{Koszulformula}
\begin{array}{rl}
2g(D_{X}Y,Z)=& Xg(Y,Z)+Yg(Z,Y)-Zg(X,Y) \\
             &-g(X,[Y,Z])+g(Y,[Z,X])+g(Z,[X,Y]).
\end{array}
\end{equation}
\end{defn}

\begin{rem}
We note that, when $g$ is \emph{non-degenerate}, $D$ is nothing but the \emph{Levi-Civita connection},
and it is uniquely determined by (\ref{Koszulformula}) (cf. \cite{ON}, Theorem 11, p.61). However,
when $g$ is \emph{degenerate}, the Koszul derivative is only determined up to a smooth section of the radical of $M$,
in the sense that, given any two Koszul derivatives $D,D'$ on $M$ and any two vector fields $X,Y\in \Gamma (TM)$,
we have $D_{X}Y-D'_{X}Y\in \Gamma ({\cal N})$.
\end{rem}

We have the following fundamental lemma of degenerate semi-Riemannian geometry:
\begin{lem}[\cite{K}, Lemma 3.1.2]
Let $(M,g)$ be a  semi-Riemannian manifold. Then
$(M,g)$ admits a Koszul derivative if and only if it is stationary.
\end{lem}
\endproof
For a later use, given an endomorphism $\sigma \in \Gamma
(\End (TM))$ of the tangent bundle $TM$, we define its Koszul derivative
by the Leibniz rule:
\begin{equation}\label{sillenzio}
(D\sigma )(Y):= D(\sigma (Y))-\sigma (DY), \qquad (Y\in \Gamma (TM)).
\end{equation}

It is easy to see that given a Koszul derivative $D$ on $M$, then
\begin{equation}\label{ninzo}
D_{X}A\in \Gamma ({\cal N}) \qquad (X\in \Gamma (TM), \, A \in \Gamma ({\cal N}) )
\end{equation}
In fact, for any $Z\in \Gamma (TM)$
we have
\begin{displaymath}
g(D_{X}A,Z)=X(g(A,Z))-g(A,D_{X}Z)=0,  \qquad (X\in \Gamma (TM), \, A \in \Gamma ({\cal N}) ).
\end{displaymath}

We have that:
\begin{lem}[\cite{K}, Lemma 3.1.4]\label{porcospino}
If the manifold $(M,g)$ is stationary then $\cal N$ is integrable.
\end{lem}
\proof
Let $A,B \in \Gamma ({\cal N})$ and let $D$ be a Koszul derivative on $M$.
Then, for any $V\in \Gamma (TM)$:
\begin{displaymath}
\begin{array}{rl}
g([A,B],V)=& g(D_{A}B,V)-g(D_{B}A,V) \\
=& A(g(B,V))-g(B,D_{A}V) - B(g(A,V))+g(A,D_{B}V) =0,
\end{array}
\end{displaymath}
so that $[A,B]\in \Gamma ({\cal N})$.
\endproof
By the Frobenius Theorem, we obtain a \emph{foliation} associated to ${\cal N}$;
we shall call this the \emph{radical foliation of $M$}.

Let $(M,g)$ be a  stationary semi-Riemannian manifold
of (constant) signature $(r,p,q)$, with $r\geq 0$. Let $E \rightarrow M$ be a
\emph{semi-Riemannian bundle} (i.e. a bundle whose
fibres are semi-Euclidean spaces of (constant) signature $(r,p,q)$); by
$\overline{E}$ (cf. (\ref{biondabella})) we shall denote the quotient
\begin{displaymath}
\overline{E}:=E\big/ {\cal N}(E)\equiv \cup _{x\in M}E_{x}\big/ {\cal N}(E_{x}),
\end{displaymath}
$E_{x}$ being the fibre of $E$ over $x\in M$. In particular,
we define the \emph{quotient tangent bundle of $M$} by
$\overline{TM}:=TM/{\cal N}(TM)$; this is endowed with
the non-degenerate metric $\overline{g}(\overline{X},\overline{Y}):=g(X,Y), \quad (X,Y\in \Gamma (TM)$)
of signature $(0,p,q)$. Let
$\overline{TM}\, ^{*} \, (=\overline{T^{*}M})$ be its dual bundle.
\begin{defn}\label{ridimentecato}
We shall call an $E$-valued 1-form $\sigma \in \Gamma
(T^{*}M\otimes E)$ {\em
radical-preserving} if, for each $x\in M$,
\begin{displaymath}
\sigma _{x}({\cal N}(T_{x}M))\subseteq {\cal N}(E_{x}).
\end{displaymath}
Denote by $\pi _{TM}:TM \rightarrow \overline{TM}$ and $\pi _{E}:E \rightarrow \overline{E}$
the natural projections. Then there exists a linear bundle map $\overline{\sigma}\in \Gamma (\overline{T^{*}M}\otimes
\overline{E})$ such that the following diagram commutes
\begin{displaymath}
\begin{array}{ccc}
TM &\stackrel{\sigma}{\longrightarrow} & E\\
\Big\downarrow\vcenter{\rlap{$\scriptstyle{\pi _{TM}}$}} & &\Big\downarrow\vcenter{\rlap{$\scriptstyle{\pi _{E}}$}} \\
\overline{TM} & \stackrel{\overline{\sigma}}{\longrightarrow} &\overline{E}
\end{array}
\end{displaymath}
if and only if $\sigma$ is radical-preserving.

We shall say that a map $\phi :M\rightarrow N$ is \emph{radical-preserving} if its differential
$d\phi \in \Gamma (T^{*}M\otimes \phi ^{-1}TN)$ is radical-preserving.
\end{defn}

We now state the fundamental theorem of singular semi-Riemannian geometry.
\begin{thm}[\cite{K}, Theorem 3.2.3]\label{fundthm}
Let $(M,g)$ be a  semi-Riemannian manifold. If
$(M,g)$ is stationary, then there exists a unique connection
$\overline{\nabla}$ on $(\overline{TM},\overline{g})$ which is torsion-free in the sense that
$\overline{T}^{\overline{\nabla}}(X,Y):=
\overline{\nabla}_{X}\overline{Y}-\overline{\nabla}_{Y}\overline{X}-\overline{[X,Y]}=0
\quad \big( X,Y\in \Gamma (TM)\big)$, and compatible with the
metric $\overline{g}$ in the sense that $\overline{\nabla} \, \overline{g}=0$; in fact
$\overline{\nabla}$ is given by:
\begin{displaymath}
\overline{\nabla}_{X}\overline{Y}:=\overline{D_{X}Y} 
\qquad \big( X\in \Gamma (TM), \, \overline{Y} \in \Gamma (\overline{TM}) \big) ,
\end{displaymath}
where $D$ is any Koszul derivative on $(M,g)$ and $Y\in \Gamma (TM)$
has $\pi _{TM}(Y)=\overline{Y}$.

Conversely, if there exists such a connection $\overline{\nabla}$, then
$(M,g)$ is stationary.
\end{thm}
\endproof

The connection $\overline{\nabla}$ is called the \emph{Koszul connection on $(M,g)$}. If $(M,g)$ is non-degenerate,
then $\overline{\nabla}$ coincides with the usual Levi-Civita connection.
Let us set $E\equiv TM$ and let $\sigma \in
\Gamma (T^{*}M \otimes TM)$ be radical-preserving. We define the Koszul connection on
$\overline{T^{*}M}\otimes \overline{TM}$ by the Leibniz rule
\begin{displaymath}
(\overline{\nabla} _{X}\overline{\sigma} )\overline{Y}:= \overline{\nabla} _{X}
(\overline{\sigma}(\overline{Y}))-\overline{\sigma}(\overline{\nabla}_{X}\overline{Y}).
\end{displaymath}
where  $X,Y \in \Gamma (TM)$ and $\overline{\nabla }$ is defined as in Theorem \ref{fundthm}.

We note that the connection $\overline{\nabla}$ is defined for $(X,\overline{Y}) \in TM\otimes \overline{TM}$, as is the
operator $\overline{\nabla}\overline{\sigma}$ defined above. It does not, in general, factor to an
operator on $\overline{TM}\otimes \overline{TM}$. However, if $\sigma =d\phi$, i.e. if
$\sigma$ is the differential of a map $\phi :M\rightarrow N$, with $\phi$ radical-preserving, we have the following
fact. Let $\phi ^{-1}(\overline{TN})\rightarrow M$ denote the pull-back of the bundle $\overline{TN} \rightarrow N$,
equivalently
\begin{displaymath}
\phi ^{-1}(\overline{TN}):=\phi ^{-1}(TN)\big/ \phi ^{-1}({\cal N}(TN)).
\end{displaymath}
\begin{lem}\label{libianconnection}
The operator $\overline{B} ^{\phi} \in \Gamma \big( \! \otimes ^{2}\overline{TM}^{*}\otimes \phi ^{-1}(\overline{TN})\big)$ defined by
\begin{displaymath}
\overline{B} ^{\phi} (\overline{X},\overline{Y})\equiv (\overline{\nabla} \, \overline{d\phi})(\overline{X}, \overline{Y}):=
(\overline{\nabla} _{X}\overline{d\phi})(\overline{Y}), \quad 
\big( \overline{X},\overline{Y} \in \Gamma (\overline{TM})\big) ,
\end{displaymath}
is well-defined, tensorial and symmetric.
\end{lem}
\endproof
We shall call the operator $\overline{B} ^{\phi}$ the \emph{(generalized) second fundamental
form of the map $\phi$}.
\section{Generalized harmonic maps and morphisms}
Let $\phi :M\rightarrow N$ be a ($C^{1}$) radical-preserving map.
We define the \emph{(generalized) differential of $\phi$} (cf. Definition \ref{ridimentecato}),
$\overline{d\phi}:\overline{TM} \rightarrow \overline{TN}$, to be the map
\begin{equation}\label{justapage}
\overline{d\phi}(\overline{X}):=\overline{d\phi (X)}, \quad \textrm{ for any } X\in \Gamma (TM).
\end{equation}
We shall define the \emph{(generalized) divergence $\overline{\diver} \, (\overline{d\phi})$} of
$\overline{d\phi}$. 
Let $\{ e_{i} \} _{i=1}^{m}$ be any basis of $TM$ such that 
${\cal N}(TM)=\spa (e_{1}, \ldots ,e_{r})$ and let
${\cal V}_{1}:=\spa (e_{r+1}, \ldots ,e_{m})$
be a \emph{screen space}, i.e. a subbundle of $TM$ such that $TM={\cal N}(TM)\oplus {\cal V}_{1}$;
we shall call such a basis a \emph{(local) radical basis for $TM$}. Then
\begin{displaymath}
\overline{\diver} \, (\overline{d\phi}) := \tr _{\overline{g}}(\overline{B} ^{\phi} ):=
\sum^{m}_{a,b=r+1} \overline{g} ^{ab}\left( \overline{\nabla} _{e_{a}}\overline{d\phi} \right) \overline{e_{b}} \, ,
\end{displaymath}
where $\overline{g}_{ab}:=\overline{g}(\overline{e_{a}},\overline{e_{b}})$.
This is well defined and does not depend on the choice of the local radical
basis $\{ e_{i} \} _{i=1}^{m}$ on $M$.

We can now define the \emph{(generalized) tension field $\overline{\tau}(\phi)$} of
a ($C^{2}$) radical-preserving map $\phi :M \rightarrow N$ between  stationary manifolds by:
\begin{displaymath}
\overline{\tau}(\phi):=\overline{\diver} \, (\overline{d\phi}).
\end{displaymath}
\begin{defn}\label{bar-harmonicOK}
We shall say that a radical-preserving map $\phi :M \rightarrow N$ between
stationary semi-Riemannian manifolds is \emph{(generalized) harmonic} if
its (generalized) tension field $\overline{\tau}(\phi )$ is identically zero.
\end{defn}
Note that this notion agrees with the usual notion of harmonicity when the
manifolds $M$ and $N$ are both non-degenerate.

If $(x^{1}, \ldots ,x^{m})$ and $(y^{1}, \ldots , y^{n})$ are \emph{radical coordinates}
(i.e. coordinates whose tangent vector fields form a radical basis) 
on $M$ and $N$ respectively (with $\rank {\cal N}(TM)=r$ and $\rank {\cal N}(TN)=\rho $), 
then, analogously to the non-degenerate case, the (generalized) tension
field of $\phi$ can be locally expressed by (cf. \cite{E-L})
\begin{equation}\label{localtension}
\overline{\tau}^{\gamma}(\phi )=\sum ^{n}_{\alpha, \beta ,\gamma =\rho +1}\sum ^{m}_{i,j,k=r+1}
\overline{g}^{ij} \left(
\phi ^{\gamma}_{ij}-{}^{M}\overline{\Gamma}^{k}_{ij}\phi ^{\gamma}_{k} +{}^{N}\overline{\Gamma}^{\gamma}_{\alpha \beta}
\phi ^{\alpha}_{i}\phi ^{\beta}_{j}\right) ,
\end{equation}
where $\phi ^{\gamma}_{k}:=\partial \phi ^{\gamma}/\partial x^{k}$, and
${}^{M}\overline{\Gamma}^{k}_{ij}\, \overline{\partial /\partial x^{k}}
:=\overline{\nabla} ^{M}_{\partial /\partial x^{i}}\, \overline{\partial /\partial x^{j}}, \,
{}^{N}\overline{\Gamma}^{\gamma}_{\alpha \beta}\, \overline{\partial /\partial y^{\gamma}}
:=\overline{\nabla} ^{N}_{\partial /\partial y^{\alpha}}\, \overline{\partial /\partial y^{\beta}}$.
In particular, if $N\equiv \RR$, then $\overline{\tau}$ reduces to what we shall call the
\emph{(generalized) Laplace-Beltrami operator} $\overline{\Delta ^{M}}$ and
the radical-preserving functions $f\in C^{\infty}(M)$ satisfying
$\overline{\Delta ^{M}}f=0$ will be called \emph{(generalized) harmonic functions}.
\begin{exam}\label{N=RR^n}
If $N=\RR ^{n}_{\rho ,\pi ,\sigma}$ then a map $\phi :(M,g)\rightarrow 
\RR ^{n}_{\rho ,\pi ,\sigma}$ is (generalized) harmonic if and only if each component
$\phi ^{\alpha}:(M,g)\rightarrow \RR , \; \alpha =1, \ldots, n$, is a (generalized) harmonic function.
\end{exam}
Now we can state the following
\begin{defn}\label{morphismOK}
We shall call a ($C^{2}$) radical-preserving map $\phi :(M,g)\rightarrow (N,h)$ between
semi-Riemannian manifolds a \emph{(generalized) harmonic morphism} if, for any (generalized) harmonic function
$f:V \subseteq N \rightarrow \RR$ on an open subset $V\subseteq N$, with $\phi ^{-1}(V)$ non-empty,
its pull-back $\phi ^{\ast}f:=f\circ \phi$ is a (generalized) harmonic function on $M$.
\end{defn}

Note that the usual definition of harmonic morphism does not make sense
for degenerate manifolds since the trace and the Laplacian
are not defined when the metric is degenerate.

Let $\phi :M \rightarrow N$ be
a radical-preserving map between two semi-Riemannian manifolds
and $\overline{d\phi} _{x}:\overline{T _{x}M} \rightarrow \overline{T _{\phi (x)}N}$ its (generalized) differential
at $x\in M$ (cf. (\ref{justapage}));
then we define the \emph{(generalized) adjoint} 
$\overline{d\phi}^{*} _{\phi (x)}:\overline{T _{\phi (x)}N} \rightarrow \overline{T _{x}M}$ of $d\phi$ as the 
adjoint of $\overline{d\phi}_{x}$, i.e. the linear map characterized by
\begin{equation}\label{bar-adjoint}
\overline{g}_{x}(\overline{d\phi}^{*}_{x}(\overline{V}), \overline{X})=
\overline{h}_{\phi (x)}(\overline{V},\overline{d\phi}_{x} (\overline{X}))=\, h_{\phi (x)}(V,d\phi _{x}(X)),
\end{equation}
for any $V\in T_{\phi (x)}N$ and $X\in T_{x}M$.

We now generalize the notion of horizontal weak conformality.
\begin{defn}\label{barHWCdefn}
We shall call a radical-preserving map $\phi :(M,g) \rightarrow (N,h)$ between two non-null
semi-Riemannian manifolds $M$ and $N$ \emph{(generalized) horizontally (weakly) conformal} (or, for
brevity, \emph{(generalized) HWC})
at $x\in M$ with square dilation $\Lambda (x)$ if
\begin{equation}\label{bar-HWC}
\overline{g}_{x}(\overline{d\phi}^{*}_{x}(\overline{V}),\overline{d\phi}^{*}_{x}(\overline{W}))=
\Lambda (x)\, \overline{h}_{\phi (x)}(\overline{V},\overline{W}), \qquad \big(
\overline{V}, \overline{W}\in \overline{T_{\phi (x)}N} \big) .
\end{equation}
In particular, if $\Lambda$ is identically equal to 1, we shall say that $\phi$ is a \emph{(generalized) Riemannian
submersion}.
\end{defn}

\begin{rem}
If both $M$ and $N$ are non-degenerate, then the above notion of (generalized) horizontal weak
conformality coincides with the better-known
one of horizontal weak conformality.
\end{rem}

Let $(M,g)$ and $(N,h)$ be stationary manifolds of signatures $\sign g=(r,p,q)$ and 
$\sign h=(\rho ,\pi ,\eta )$, respectively,
and let $\phi :M\rightarrow N$ be a radical-preserving map (i.e. a map whose differential
$d\phi _{x}$ is radical-preserving for each $x\in M$).
As usual, for any $x\in M$, set
${\cal V}_{x}:=\ker d\phi _{x}$ and ${\cal H}_{x}:={\cal V}_{x}^{\perp}$. We shall also set:
\begin{displaymath}
\overline{\cal V}_{x}:={\cal V}_{x}/({\cal N}(T_{x}M)\cap {\cal V}_{x}), \quad
\overline{\cal H}_{x}:={\cal H}_{x}/{\cal N}(T_{x}M),
\end{displaymath}
having noticed that, by equation (\ref{nig}), ${\cal N}(T_{x}M)\subseteq {\cal H}_{x}$. We have the following
\begin{lem}\label{imagedphistarlem}
Let $\phi :M\rightarrow N$ be a radical-preserving map.
Then, at any $x\in M$, the following identity holds:
\begin{equation}\label{imagedphistar2}
\image \overline{d\phi }_{x}^{*} =\overline{\cal H}_{x}\, .
\end{equation}
\end{lem}
\proof
For any $x\in M$, it is easy to see that $\ker \overline{d\phi}_{x}=\overline{\ker d\phi}_{x}$.
Then we have
\begin{displaymath}
\image \overline{d\phi}^{*}_{x}=(\ker \overline{d\phi}_{x})^{\perp _{\overline{g}}}=
(\overline{\cal V}_{x})^{\perp _{\overline{g}}}=\overline{\cal H}_{x}\, ,
\end{displaymath}
the last equality following by Proposition \ref{bigdiscovery}.
\endproof
We have the following special sort of generalized HWC maps
\begin{lem}\label{HsubseteqVlem}
Let $\phi :M \rightarrow N$ be a radical-preserving map between non-null
semi-Riemannian manifolds.
Then $\phi$ is (generalized) HWC at $x\in M$ with square dilation $\Lambda (x)=0$ if and only if
\begin{equation}\label{HsubseteqV}
\overline{\cal H}_{x}\subseteq \overline{\cal V}_{x} \, ,
\end{equation}
i.e. if and only if $\overline{\cal H}_{x}$ is null.
\end{lem}
\proof
By Definition \ref{barHWCdefn}, $\phi$ is (generalized) HWC with square dilation $\Lambda (x)=0$ if and only if
\begin{displaymath}
\overline{g}_{x}(\overline{d\phi }_{x}^{*}(\overline{V}),\overline{d\phi }_{x}^{*}(\overline{W}))=0
\qquad (V,W \in T_{\phi (x)}N).
\end{displaymath}
By equation (\ref{imagedphistar2}), this holds if and only if
$\overline{\cal H}_{x}$ is null, equivalently, bearing in mind (\ref{imagedphistar2}), equation (\ref{HsubseteqV}) holds.
\endproof

We have the following characterization which generalizes a better-known characterization of HWC maps (cf. \cite{B-W}).
\begin{prop}\label{doveseiandata}
A radical-preserving map $\phi :M \rightarrow N$ between non-null semi-Riemannian manifolds 
$M$ and $N$ is (generalized) HWC at $x\in M$ with
square dilation $\Lambda (x)$ if and only if
\begin{equation}\label{dphidphistar=Lambda}
\overline{d\phi}_{x} \circ \overline{d\phi}^{*}_{x}=\Lambda (x) \, \mathbf{1}_{\overline{T_{\phi (x)}N}}\, .
\end{equation}
\end{prop}
\proof
From the characterization (\ref{bar-adjoint}) of the adjoint map $\overline{d\phi}^{*}_{x}$, we have
\begin{equation}
\overline{g}_{x}(\overline{d\phi}^{*}_{x}(\overline{V}),\overline{d\phi}^{*}_{x}(\overline{W}))=
\overline{h}_{\phi (x)}(\overline{V},\overline{d\phi}_{x} \circ \overline{d\phi}^{*}_{x}(\overline{W})),
\qquad (V,W\in T_{\phi (x)}N).
\end{equation}
Comparing with equation (\ref{bar-HWC}), gives the statement.
\endproof

\begin{prop}\label{james}
If $\phi$ is (generalized) HWC, then $\overline{\cal H}_{x} \subseteq
 \overline{\cal V}_{x}$ if and only if one of the
following holds:
\begin{enumerate}
\item [(i)] $\ker \overline{d\phi}_{x} \equiv \overline{T_{x}M}$;
\item [(ii)] $\ker \overline{d\phi}_{x}\varsubsetneqq \overline{T_{x}M}$ is degenerate.
\end{enumerate}
\end{prop}
\proof
In fact, if $\overline{\cal H}_{x} \subseteq
 \overline{\cal V}_{x}$ and \emph{(i)} does not hold, then $\overline{\cal H} _{x}\neq \{ \mathbf{0} \}$, so that there exists a vector
 $\overline{0} \neq \overline{X}\in 
\overline{\cal H} _{x}$ and, for such a vector, $\overline{g}(\overline{X}, \overline{V})=0$ for any $\overline{V}\in 
\overline{\cal V} _{x}$, so that \emph{(ii)} holds.
Conversely if $\overline{\cal V}_{x}:=\ker \overline{d\phi}_{x}\equiv \overline{T_{x}M}$ then clearly
$\overline{\cal H}_{x} \subseteq \overline{\cal V}_{x}$. If, on the other hand, $\ker \overline{d\phi}_{x}$ is
degenerate, then since $\phi$ is (generalized) HWC, we get $\Lambda (x)=0$; in fact,
$\ker \overline{d\phi}_{x}$ is degenerate if and only if $\overline{\cal H}_{x}$ is degenerate
if and only if $\overline{\cal V}_{x}\cap \overline{\cal H}_{x}\neq \{ \mathbf{0} \}$, so that
there exists a non-zero vector $\overline{V}\in \overline{T_{\phi (X)}N}$ such that
\begin{displaymath}
\overline{0}\neq \overline{d\phi}^{*}_{x}(\overline{V})\in \ker \overline{d\phi}_{x} \cap \image
\overline{d\phi}^{*}_{x}.
\end{displaymath}
Combining this with the (generalized) HWC condition gives:
\begin{displaymath}
0=\overline{g}(\overline{d\phi}^{*}_{x}(\overline{V}),\overline{d\phi}^{*}_{x}(\overline{W}))=
\Lambda (x)\, \overline{h} (\overline{V},\overline{W}) \quad \textrm{for any }\overline{W}\in
\overline{T_{\phi (X)}N},
\end{displaymath}
and, as $\overline{h}$ is non-degenerate, we must have $\Lambda (x)=0$.
Then, from Lemma
\ref{HsubseteqVlem}, $\overline{\cal H}_{x} \subseteq \overline{\cal V}_{x}$, and this gives the claim.

When the metrics $g$ and $h$ are both degenerate, the case (i) splits into the two subcases:
\begin{enumerate}
\item [(i')] $\ker d\phi _{x} \equiv T_{x}M$, i.e. $d\phi _{x} = 0$ or
\item [(i'')] $d\phi _{x} \neq 0$ and $\overline{d\phi}_{x}= 0$, i.e. 
$\{ \mathbf{0} \} \neq \image (d\phi _{x}) \subseteq {\cal N}(T_{\phi (x)}N)$.
\end{enumerate}
\endproof

In the case when the square dilation is non-zero, we have the following characterization:
\begin{prop}\label{Lambdanonzero}
A map $\phi :(M,g)\rightarrow (N,h)$ between non-null semi-Riemannian manifolds $M$ and $N$
is (generalized) HWC at $x\in M$ with square dilation 
$\Lambda (x) \neq 0$ if and only if
\begin{equation}\label{DDD}
h_{\phi (x)}(d\phi _{x}(X),d\phi _{x}(Y))=\Lambda (x)\, g_{x}(X,Y), \qquad (X,Y \in {\cal H}_{x}).
\end{equation}
\end{prop}
\proof
Suppose that $\phi$ is (generalized) HWC; then by Lemma \ref{imagedphistarlem} we have 
$\image (\overline{d\phi}^{*}_{x})=\overline{\cal H}_{x}$, so that for any
$\overline{X}, \overline{Y} \in \overline{\cal H}_{x}$ there exist vectors
$\overline{V}$ and $\overline{W} \in \overline{T_{\phi (x)}N}$ such that
\begin{equation}\label{XV}
\overline{d\phi}^{*}_{x}(\overline{V})=\overline{X} \textrm{ and }
\overline{d\phi}^{*}_{x}(\overline{W})=\overline{Y}.
\end{equation}
Applying the operator $\overline{d\phi}_{x}$ to both sides
of the identities (\ref{XV}),
and using equation (\ref{dphidphistar=Lambda}), since $\Lambda (x)\neq 0$ we obtain
\begin{displaymath}
\overline{V}=(\Lambda (x))^{-1}\overline{d\phi (X)} \textrm{ and }
\overline{W}=(\Lambda (x))^{-1}\overline{d\phi (Y)};
\end{displaymath}
on substituting these into the definition of (generalized) HWC, we obtain the statement.
The converse is similar.
\endproof
We thus obtain the following characterization for a (generalized) HWC map (cf. \cite{B-W} for
the non-degenerate case):
\begin{thm}\label{bigguy}
Let $\phi :M \rightarrow N$ be a radical-preserving map between non-null semi-Riemannian manifolds $(M,g)$
and $(N,h)$. Then $\phi$ is (generalized) HWC at $x\in M$, with square
dilation $\Lambda (x)$, if and only if precisely one of the following possibilities holds:
\begin{enumerate}
\item [(a)] $d\phi _{x}=0$ (so $\Lambda (x)=0$);
\item [(b)] $d\phi _{x}\neq 0$ and $\image (d\phi _{x}) \subseteq {\cal N}(T_{\phi (x)}N)$
(so $\Lambda (x)=0$ and $\overline{d\phi}_{x}=\overline{0}$);
\item [(c)] $\overline{\cal V}_{x}\varsubsetneqq \overline{T_{x}M}$ is degenerate and $\overline{\cal H}_{x} \subseteq \overline{\cal V}_{x}$
           (equivalently $\overline{\cal H}_{x}$ is non-zero and null): then $\Lambda (x)=0$ but $\overline{d\phi }_{x} \neq 0$;
\item [(d)] $\Lambda (x) \neq 0$ and
\begin{displaymath}
h_{\phi (x)}(d\phi _{x}(X),d\phi _{x}(Y))=\Lambda (x)\, g_{x}(X,Y) \qquad (X,Y \in {\cal H}_{x}).
\end{displaymath}
\end{enumerate}
\end{thm}
\proof
Let $x\in M$ and suppose that $\phi$ is (generalized) HWC at $x$. 
If $\Lambda (x)=0$ then by Lemma \ref{HsubseteqVlem}
we have $\overline{\cal H}_{x}\subseteq \overline{\cal V}_{x}$, so by Proposition \ref{james}, either
(i) $\ker\overline{d\phi _{x}}\equiv \overline{T_{x}M}$ or (ii) $\ker \overline{d\phi _{x}}
\varsubsetneqq \overline{T_{x}M}$ is
degenerate; the possibility (i) holds if and only if either \emph{(a)} or \emph{(b)} holds, whereas (ii) occurs
if and only if case \emph{(c)} holds.
Otherwise $\Lambda (x)\neq 0$, so that by Proposition \ref{Lambdanonzero} we obtain case
\emph{(d)}.

Conversely, if \emph{(a)}, \emph{(b)} or \emph{(c)} holds, then clearly $\phi$ is (generalized) HWC at $x$ with $\Lambda (x)=0$.
If \emph{(d)} holds then, by Proposition \ref{Lambdanonzero}, $\phi$ is (generalized) HWC at $x$ with square
dilation $\Lambda (x)\neq 0$.
\endproof

Thus the possibility of the metrics $g$ and $h$ being degenerate has given rise to another type of
point (case \emph{(b) }) which is not possible when $M$ and $N$ are non-degenerate.

We have the following characterization of (generalized) horizontal weak conformality
whose proof is similar to its (non-degenerate semi-)Riemannian analogue (cf. \cite{B-W}):
\begin{lem}\label{var-barHWClemma}
A radical-preserving map $\phi :(M,g) \rightarrow (N,h)$ between stationary manifolds is (generalized) horizontally weakly
conformal at a point $x\in M$ with square dilation $\Lambda (x)$ if and only if, in radical coordinates
$\{ x^{j} \} _{j=1}^{m}$ in a neighbourhood of $x\in M$ and $\{ y^{\alpha} \}_{\alpha =1}^{n}$ around
$\phi (x) \in N$, we have
\begin{equation}\label{var-barHWC}
\phi ^{\alpha}_{i}\phi ^{\beta}_{j}\overline{g}^{ij}=\Lambda (x)\, \overline{h}^{\alpha \beta},
\end{equation}
where $r+1 \leq i,j \leq m$, $\rho +1 \leq \alpha ,\beta \leq n$ and $\phi ^{\gamma}_{k}:=
\partial \phi ^{\gamma}/\partial x^{k}$. Moreover, setting $\overline{\grad}\phi ^{\alpha}:=
\overline{g}^{ij}\phi ^{\alpha}_{i}\overline{\partial /\partial x^{j}}$, equation (\ref{var-barHWC})
above reads:
\begin{equation}\label{var-barHWC2}
\overline{g}(\overline{\grad}\phi ^{\alpha},\overline{\grad}\phi ^{\beta})=\Lambda (x)\,
\overline{h}^{\alpha \beta}.
\end{equation}
\end{lem}

\section{A Fuglede--Ishihara-type characterization of (generalized) harmonic morphisms}

\subsection{Preliminaries}
Recall (see \cite{Mol}) that 
(i) a foliation $\cal F$ on a manifold $M$ is said to be \emph{simple} if its leaves are the (connected) fibres of
a smooth submersion defined on $M$;
(ii) the \emph{leaf space} of a foliation $\cal F$ is the topological space $M/{\cal F}$, equipped with the
quotient topology.

We note that this space, in general, is not Hausdorff. However, the following holds.
\begin{prop}\label{phwaaarh}
A foliation $\cal F$ on $M$ is simple if and only if its leaf space $M/{\cal F}$ can be given the
structure of a Hausdorff (smooth) manifold such that the natural projection $M\rightarrow M/{\cal F}$ is a
smooth submersion. Furthermore, if such a smooth structure exists, then it is unique.
\end{prop}
\endproof
Since each point $x\in M$ has a neighbourhood $W\subseteq M$ with ${\cal F}|_{W}$ simple,
${\cal F}$ is always simple \emph{locally}. Hence, as all the considerations in this
section will be \emph{local}, by replacing the manifold $M$ by a suitable open subset $W$ if necessary, \emph{we shall assume
that any foliation $\cal F$ on $M$ is simple}. We make the same assumption for $N$.

We recall (cf. Lemma \ref{porcospino}) that, if a manifold $M$ is stationary, then
its radical distribution ${\cal N}(TM)$ is integrable. Let ${\cal F}^{M}$ be the
radical foliation of $M$ (i.e. the foliation whose leaves are
tangent to ${\cal N}(TM)$); set $\overline{M}:=M/{\cal F}^{M}$, the leaf space of ${\cal N}(TM)$, and denote by
$\pi _{M}:M\rightarrow \overline{M}$ the natural projection; by
Proposition \ref{phwaaarh}, $\overline{M}$ is a smooth manifold. Elements of $\overline{M}$ will
be denoted by $[x]_{{\cal F}^{M}}:=\pi _{M}(x)$, where $x\in M$. Then, any radical-preserving map $\phi :M\rightarrow N$ between
stationary manifolds factors to a map $\overline{\phi}: \overline{M} \rightarrow \overline{N}$
in the sense that the following diagram commutes
\begin{displaymath}
\begin{array}{ccc}
M &\stackrel{\phi}{\longrightarrow} & N\\
\Big\downarrow\vcenter{\rlap{$\scriptstyle{\pi _{M}}$}} & &\Big\downarrow\vcenter{\rlap{$\scriptstyle{\pi _{N}}$}} \\
\overline{M} & \stackrel{\overline{\phi}}{\longrightarrow} &\overline{N}
\end{array} \qquad .
\end{displaymath}
Thus $\overline{\phi}([x]_{{\cal F}^{M}}):=[\phi (x)]_{{\cal F}^{N}}$. For any $[x]\in \overline{M}$,
the map $\overline{\phi}$ naturally induces a linear operator
$(d\overline{\phi})_{[x]}:T_{[x]}\overline{M} \rightarrow T_{\overline{\phi}([x])}\overline{N}$.

For each $x\in M$ define a following map
\begin{displaymath}
\Psi ^{M}_{x}: \overline{T_{x}M} \rightarrow T_{\pi _{M}(x)}\overline{M}, \qquad
\overline{X} \mapsto (d\pi _{M})_{x}(X),
\end{displaymath}
where $X\in T_{x}M$ is such that $\pi _{TM}(X)=\overline{X})$.
It is easy to see that $\Psi ^{M}_{x}$ is a well-defined isomorphism, and that the following holds
\begin{lem}
Let $\phi :M\rightarrow N$ be a radical-preserving map between stationary manifolds; then, for any $x\in M$,
\begin{equation}\label{stellina}
\Psi ^{N}_{\phi (x)}\circ \overline{d\phi}_{x}=(d\overline{\phi})_{[x]}\circ \Psi ^{M}_{x}\, ,
\end{equation}
equivalently, the following diagram commutes:
\begin{displaymath}
\begin{array}{ccc}
\overline{T_{x}M} &\stackrel{\overline{d\phi}_{x}}{\longrightarrow} & \overline{T_{\phi(x)}N}\\
\Big\downarrow\vcenter{\rlap{$\scriptstyle{\Psi ^{M}_{x}}$}} & &\Big\downarrow\vcenter{\rlap{$\scriptstyle{\Psi ^{N}_{\phi (x)}}$}} \\
T_{[x]}\overline{M} & \stackrel{(d\overline{\phi})_{[x]}}{\longrightarrow} &T_{\overline{\phi}([x])}\overline{N}
\end{array}\quad .
\end{displaymath}
In particular, as the maps $\Psi ^{M}_{x}$ and $\Psi ^{N}_{\phi (x)}$ are isomorphisms, we can identify
$(\overline{d\phi})_{x}$ with $(d\overline{\phi})_{[x]}$.
\end{lem}
\endproof
\subsection{Horizontal weak conformality of $\overline{\phi}$}
Let $(M,g)$ be a non-null stationary semi-Riemannian manifold. Then we can endow $\overline{M}$ with the induced metric 
$g^{\overline{M}}$ defined by:
\begin{displaymath}
g^{\overline{M}}:=((\Psi ^{M})^{-1})^{*}\overline{g},
\end{displaymath}
where $\overline{g}$ is defined by
\begin{displaymath}
\overline{g}(\overline{X},\overline{Y}):=g(X,Y) \qquad (X,Y\in \Gamma (TM)).
\end{displaymath}
Note that the metric $g^{\overline{M}}$ is non-degenerate.

The adjoint of $d\overline{\phi}_{[x]}:T_{[x]}\overline{M}\rightarrow T_{\overline{\phi}([x])}\overline{N}$
is the (unique) linear map $(d\overline{\phi})^{*}_{[x]}:T_{\overline{\phi}([x])}\overline{N}
\rightarrow T_{[x]}\overline{M}$ characterized as usual by
\begin{equation}
g^{\overline{M}}_{[x]}((d\overline{\phi})^{*}_{[x]}(\widetilde{V}),\widetilde{X})=
h^{\overline{N}}_{\overline{\phi}([x])}(\widetilde{V},d\overline{\phi}_{[x]}(\widetilde{X})), \quad
\quad (\widetilde{X}\in \Gamma (T_{[x]}\, \overline{M}), \widetilde{V}\in \Gamma (T_{\overline{\phi}([x])}\, \overline{N})).
\end{equation}
Setting $\widetilde{X}=\Psi ^{M}(\overline{X})$ and $\widetilde{V}=\Psi ^{M}(\overline{V})$ for some
$\overline{X}\in \Gamma (\overline{TM}), \overline{V}\in \Gamma (\overline{TN})$ and using
equation (\ref{stellina}) we obtain:
\begin{equation}\label{calzelunghe}
(d\overline{\phi})^{*}\circ \Psi ^{N}=\Psi ^{M}\circ \overline{d\phi}^{*}.
\end{equation}
Now we can state the
\begin{prop}\label{LoriHWC}
Let $\phi :(M,g)\rightarrow (N,h)$ be a radical-preserving map between stationary manifolds. Then
$\phi$ is (generalized) HWC if and only if $\overline{\phi}$ is HWC.
\end{prop}
\proof
The map $\overline{\phi}$ is HWC with square dilation $\Lambda$ if and only if:
\begin{equation}\label{freddoeh}
g^{\overline{M}}\big( (d\overline{\phi})^{*}(\widetilde{V}),(d\overline{\phi})^{*}(\widetilde{W})\big) =\Lambda
\, h^{\overline{N}}(\widetilde{V},\widetilde{W}), \quad
\quad (\widetilde{V},\widetilde{W}\in \Gamma (T\, \overline{N})).
\end{equation}
Let $\overline{V},\overline{W}\in \Gamma (\overline{TN})$ be such that:
\begin{equation}\label{substityyy}
\widetilde{V}=\Psi ^{N}(\overline{V}), \quad \widetilde{W}=\Psi ^{N}(\overline{W});
\end{equation}
then, on using substitutions (\ref{substityyy}), equation (\ref{calzelunghe}) and the definition of $g^{\overline{M}}$, 
we see that (\ref{freddoeh}) is equivalent to $\phi$ being (generalized) HWC.
\endproof

\subsection{On harmonicity of $\overline{\phi}$}
Let $(M,g)$ and $(N,h)$ be two stationary manifolds of dimension $m$ and $n$ respectively, whose radical
distributions ${\cal N}(TM)$ and ${\cal N}(TN)$ have ranks $r$ and $\rho$ respectively.
Then the quotient manifolds $(\overline{M},g^{\overline{M}})$ and $(\overline{N},h^{\overline{N}})$
are $(m-r)$- and $(n-\rho)$-dimensional non-degenerate semi-Riemannian manifolds, thus 
they admit uniquely determined Levi-Civita connections $\nabla ^{\overline{M}}$ and
$\nabla ^{\overline{N}}$, respectively.

As $\overline{M}$ and $\overline{N}$ are non-degenerate, we have the usual notion of \emph{tension field} $\tau$, for a
map $\overline{\phi}:\overline{M} \rightarrow \overline{N}$:
\begin{equation}\label{tensionphibar}
\tau (\overline{\phi}):=\tr _{g^{\overline{M}}}(\nabla d\, \overline{\phi}),
\end{equation}
where $\nabla$ is the connection on the bundle $(T \, \overline{M})^{*}\otimes (\overline{\phi}) ^{-1}
(T \, \overline{N})$ induced from $\nabla ^{\overline{M}}$ and $\nabla ^{\overline{N}}$.
Then $\overline{\phi}$ is harmonic if and only if $\tau (\overline{\phi})=0$. 
Endow   $(M,g)$ (resp. $(N,h)$)
with (local) radical coordinates $(x^{1}, \ldots ,x^{r}, x^{r+1}, \ldots ,x^{m})$ (resp.
$(y^{1}, \ldots ,y^{\rho}, y^{\rho +1}, \ldots , y^{n})$); then
$\overline{M}$ (resp. $\overline{N}$) has the same coordinates as $M$ with the first $r$ (resp. $\rho$) coordinates
omitted.
In these coordinates, (\ref{tensionphibar}) reads:
\begin{displaymath}
\tau ^{\gamma}(\overline{\phi})=\sum ^{n}_{\alpha ,\beta ,\gamma =\rho +1}\sum _{i,j,k=r+1}^{m}
(g^{\overline{M}})^{ij}(\overline{\phi}^{\gamma}_{ij}-{}^{\overline{M}}\Gamma ^{k}_{ij}\overline{\phi}^{\gamma}_{k}+
{}^{\overline{N}}\Gamma ^{\gamma}_{\alpha \beta}\overline{\phi}^{\alpha}_{i}\overline{\phi}^{\beta}_{j}),
\end{displaymath}
where ${}^{\overline{M}}\Gamma ^{k}_{ij}\partial /\partial x^{k}:= \nabla ^{\overline{M}}_{\partial /\partial x^{i}}
\partial /\partial x^{j}$ and ${}^{\overline{N}}\Gamma ^{\gamma}_{\alpha \beta}\partial /\partial y^{\gamma}:=
\nabla ^{\overline{N}}_{\partial /\partial y^{\alpha}} \partial /\partial y^{\beta}$.
Since the coordinates are radical, we have $\overline{\phi}^{\gamma}=\phi ^{\gamma}$ (for $\gamma =\rho +1, \ldots ,n$), 
and the Christoffel symbols
${}^{\overline{M}}\Gamma ^{k}_{ij}$ and ${}^{\overline{N}}\Gamma ^{\gamma}_{\alpha \beta}$ agree
with the symbols ${}^{M}\overline{\Gamma }^{k}_{ij}$ and ${}^{N}
\overline{\Gamma }^{\gamma}_{\alpha \beta}$ of formula (\ref{localtension}) (for
$r+1 \leq i,j,k \leq m$ and $\rho +1 \leq \alpha ,\beta ,\gamma \leq n$); hence, we have:
\begin{prop}\label{cuginoLori}
Let $\phi :M\rightarrow N$ be a radical-preserving map between stationary manifolds. Then,
on identifying $\overline{T_{y}N}$ with $T_{\overline{y}}\overline{N} \, (\overline{y}:=\pi _{N}(y))$,
$\overline{\tau}(\phi)_{x}\in \overline{T_{\phi (x)}N}$ can be identified with $\tau (\overline{\phi})_{\overline{x}}\in
T_{\overline{\phi}(\overline{x})}\overline{N}$;
in particular, $\overline{\phi}$ is harmonic if and only if $\phi$ is (generalized) harmonic.
\end{prop}
\endproof
\subsection{Main characterization of (generalized) harmonic morphisms and examples}

Now we state the Fuglede--Ishihara-type characterization for (generalized) harmonic morphisms:

\begin{thm}\label{bigguy2}
Let $\phi :M\rightarrow N$ be a radical-preserving map between stationary manifolds. Then $\phi$ is a
(generalized) harmonic morphism if and  only if it is (generalized) harmonic and (generalized) HWC.
\end{thm}
\proof
Any (generalized) harmonic function $f:U\subseteq N \rightarrow \RR$ is, by definition, radical-preserving, and so factors to a 
smooth function
$\overline{f}: \pi _{N}(U)\subseteq \overline{N} \rightarrow \RR$, with $f=\overline{f} \circ \pi _{N}$; this
function $\overline{f}$ is harmonic, by Proposition \ref{cuginoLori}. Conversely, if
$\overline{f}:V\subseteq \overline{N}\rightarrow \RR$ is harmonic, then $f:=\overline{f}\circ \pi _{N}$
is (generalized) harmonic. Hence,
the map $\phi$ is a (generalized) harmonic morphism if and only if $\overline{\phi}:\overline{M}
\rightarrow \overline{N}$ is a harmonic morphism. By Fuglede's Theorem (cf. \cite{Fu96}, Theorem 3)
this is equivalent to $\overline{\phi}$ being harmonic and
HWC, then the claim follows from Propositions \ref{cuginoLori} and \ref{LoriHWC}.
\endproof

Now we give few examples of (generalized) harmonic morphisms.

\begin{exam}
Let $\phi$ be a ($C^{2}$) map
\begin{displaymath}
\phi :\RR ^{3}_{1,1,1}\backslash \{ x^{2}=x^{3} \}\rightarrow \RR ,\quad
(x^{1},x^{2},x^{3})                     \mapsto \phi (x^{1},x^{2},x^{3}).
\end{displaymath}
Clearly ${\cal N}(\RR ^{3}_{1,1,1})=\spa (\partial /\partial x^{1})$ and ${\cal N}(\RR )=\{
\mathbf{0} \} $. Moreover we have
$d\phi \left( \partial /\partial x^{1} \right) = \partial \phi /\partial x^{1}$,
so $\phi$ is radical-preserving if and only if $\partial \phi /\partial x^{1}=0$.
We notice that the coordinates $(x^{1},x^{2},x^{3})$ are radical. Identifying the vector fields
$\partial /\partial x^{2}$ and $\partial /\partial x^{3} \in \Gamma (T\RR ^{3}_{1,1,1})$ and $\partial /\partial t \in
\Gamma (T\RR )$
with their natural projections in $\overline{T\RR ^{3}_{1,1,1}}$ and $\overline{T\RR}$ respectively,
a simple calculation gives the following expression for $\overline{d\phi} ^{*}$:
\begin{displaymath}
\overline{d\phi} ^{*}\left( \frac{\partial}{\partial t} \right) = -\frac{\partial \phi}{\partial x^{2}}
\frac{\partial }{\partial x^{2}}+\frac{\partial \phi}{\partial x^{3}}
\frac{\partial }{\partial x^{3}},
\end{displaymath}
from which we get:
\begin{equation}\label{lambdaexam}
\Big\langle \overline{d\phi} ^{*}\left( \frac{\partial}{\partial t} \right) ,
\overline{d\phi} ^{*}\left( \frac{\partial}{\partial t} \right) \Big\rangle _{\overline{T\RR ^{3}_{1,1,1}}}=
-\left( \frac{\partial \phi}{\partial x^{2}} \right) ^{2}+
\left( \frac{\partial \phi}{\partial x^{3}} \right) ^{2}=:\Lambda.
\end{equation}
So $\phi$ is a (generalized) HWC with square dilation $\Lambda$.

Moreover $\phi$ is (generalized) harmonic if and only if
\begin{displaymath}
\frac{\partial ^{2}\phi}{(\partial x^{2})^{2}} -
\frac{\partial ^{2}\phi}{(\partial x^{3})^{2}}=0
\end{displaymath}
i.e. if and only if $\phi$ is of the form $\phi (x^{1},x^{2},x^{3})=\mu (x^{2}+x^{3})+\nu (x^{2}-x^{3})$, where
$\mu ,\nu \in C^{2}(\RR )$. By Theorem \ref{bigguy2}, $\phi$ is a (generalized) harmonic morphism.
\end{exam}
\begin{exam}[An anti-orthogonal multiplication]
Identify $\RR ^{3}_{1,1,1}$ with the (associative) algebra 
\begin{displaymath}
\{ x=\epsilon x^{1}+\eta x^{2}+jx^{3}, \quad
(x^{1},x^{2},x^{3})\in \RR ^{3} \} ,
\end{displaymath}
where $\epsilon ,\eta$ and $j$ satisfy the following relations:
\begin{displaymath}
\epsilon ^{2}=\epsilon \eta =\eta \epsilon =\epsilon j=j \epsilon =0, \quad
j^{2}=\eta ^{2}=\eta , \quad \eta j =j\eta =j.
\end{displaymath}
Given two elements $x,y\in \RR ^{3}_{1,1,1}$ we can define their product
\begin{displaymath}
\theta :\RR ^{3}_{1,1,1}\times \RR ^{3}_{1,1,1} \rightarrow \RR ^{2}_{0,1,1}\subseteq \RR ^{3}_{1,1,1},
\qquad \theta (x,y)=x\cdot y,
\end{displaymath}
as follows:
\begin{displaymath}
\begin{array}{rl}
\theta (x,y) &= x\cdot y \\
&=(\epsilon x^{1}+\eta x^{2}+jx^{3})(\epsilon y^{1}+\eta y^{2}+jy^{3}) \\
&=\epsilon \cdot 0+\eta (x^{2}y^{2}+x^{3}y^{3}) +j(x^{2}y^{3}+x^{3}y^{2}).
\end{array}
\end{displaymath}
For any $x\in \RR ^{3}_{1,1,1}$ we define the \emph{square norm $\Arrowvert x \Arrowvert _{1,1,1}^{2}$}
(induced from the metric on $\RR ^{3}_{1,1,1}$) by:
\begin{displaymath}
\Arrowvert x \Arrowvert _{1,1,1}^{2}:=-(x^{2})^{2}+(x^{3})^{2}.
\end{displaymath}
Then $\Arrowvert \theta (x,y) \Arrowvert _{1,1,1}^{2} =-\Arrowvert x \Arrowvert _{1,1,1}^{2}
\cdot \Arrowvert y \Arrowvert _{1,1,1}^{2}$, so $\theta$ is an \emph{anti-orthogonal multiplication}.    

Take standard coordinates $(x^{1},x^{2},x^{3},y^{1},y^{2},y^{3})$ in $\RR ^{3}_{1,1,1}\times \RR ^{3}_{1,1,1}$, and
$(z^{1},z^{2},z^{3})$ in $\RR ^{3}_{1,1,1}$. They are radical coordinates. It is easy to see that:
\begin{displaymath}
{\cal N}(\RR ^{3}_{1,1,1}\times \RR ^{3}_{1,1,1}):=\spa \left(
\frac{\partial}{\partial x^{1}},\frac{\partial}{\partial y^{1}} \right)
\end{displaymath}
and
\begin{displaymath}
{\cal N}(\RR ^{3}_{1,1,1}):=\spa \left( \frac{\partial}{\partial z^{1}} \right) .
\end{displaymath}
Moreover
\begin{displaymath}
d\theta = (0,x^{2}dy^{2}+y^{2}dx^{2}+x^{3}dy^{3}+y^{3}dx^{3},x^{2}dy^{3}+y^{3}dx^{2}+x^{3}dy^{2}+y^{2}dx^{3}),
\end{displaymath}
so that $\theta$ is radical-preserving.

The components $\theta ^{i}, \quad i=1,2,3$ of $\theta$ are easily seen to be (generalized) harmonic, so that $\theta$ is
(generalized) harmonic.

In order to check the (generalized) horizontal weak conformality, we make use of Lemma \ref{var-barHWClemma}. 
So, in this case, $\theta$ is (generalized) HWC since
\begin{displaymath}
\overline{g}^{ij}(\theta ^{\alpha}_{i})(\theta ^{\beta}_{j})=\Lambda \overline{h}^{\alpha \beta},
\end{displaymath}
where $(\overline{g}^{ij})=\diag (-1,1,-1,1)$ and $(\overline{h}^{\alpha \beta}) =
\diag (-1,1)$ and $\Lambda = -\big( -(y^{2})^{2}+(y^{3})^{2} -(x^{2})^{2}+(x^{3})^{2}\big) =
-(\Arrowvert x \Arrowvert _{1,1,1}^{2} +\Arrowvert y \Arrowvert _{1,1,1}^{2})$.
Finally, applying Theorem \ref{bigguy2}, we see that $\theta$ is a \emph{(generalized) harmonic morphism}.
\end{exam}

\begin{exam}[Radial projection]\label{radialproj}
Let $\RR ^{3}_{1,1,1}$ be $\RR ^{3}$ endowed with
the degenerate metric $g=-(dx^{2})^{2}+(dx^{3})^{2}$, where $(x^{1}, x^{2},x^{3})$ are the canonical
(and so radical) coordinates on $\RR ^{3}$. We set:
\begin{displaymath}
(\RR ^{3}_{1,1,1})^{+}:=(\RR ^{3} \backslash \{-(x^{2})^{2}+(x^{3})^{2}\leq 0 \}, g).
\end{displaymath}
We define the \emph{degenerate 2-pseudo-sphere $S^{2}_{1,1,1}$} as the manifold:
\begin{displaymath}
S^{2}_{1,1,1}:= \{ x\in \RR ^{3}: -(x^{2})^{2}+(x^{3})^{2}=1 \} ,
\end{displaymath}
endowed with the induced metric $h:=i^{*}g$, where $i:S^{2}_{1,1,1}\hookrightarrow \RR ^{3}_{1,1,1}$
is the natural inclusion. We can then define the following map:
\begin{displaymath}
\phi :(\RR ^{3}_{1,1,1})^{+} \rightarrow S^{2}_{1,1,1} \subseteq \RR ^{3}_{1,1,1}, \quad 
x  \mapsto x/\Arrowvert x \Arrowvert ,
\end{displaymath}
where $\Arrowvert x \Arrowvert :=\sqrt{-(x^{2})^{2}+(x^{3})^{2} }$ is the norm with respect
to the metric of $\RR ^{3}_{1,1,1}$.
As $\dim \overline{T_{x}S^{2}_{1,1,1}}=1$, $\phi$ is automatically (generalized) HWC.
Set $\phi ^{\alpha}_{i}:=\partial \phi ^{\alpha}/\partial x^{i} \;
(\alpha =1,2,3 \textrm{ and } i=1,2)$. From Lemma \ref{var-barHWClemma},
by parametrizing the upper half of $S^{2}_{1,1,1}$ by
$X=X(t,u):=(t, \sinh u,\cosh u )\subseteq \RR ^{3}_{1,1,1}$, we find that, for $x^{2}\neq 0$,

\begin{displaymath}
\Lambda (x)=(\phi ^{2}_{2})^{2}-(\phi ^{2}_{3})^{2}=
\frac{1}{(x^{3})^{2}} \left( 1-\left( \frac{x^{2}}{\Arrowvert x \Arrowvert} \right) ^{2} \right) ^{\! 2}-
\frac{1}{(x^{2})^{2}} \left( 1-\left( \frac{x^{3}}{\Arrowvert x \Arrowvert} \right) ^{2} \right) ^{\! 2}.
\end{displaymath}
and
\begin{displaymath}
\ker d\phi _{x}=\spa \left( \frac{x^{1}(x^{3}-\gamma x^{2})}{\Arrowvert x \Arrowvert ^{2}} \, 
\frac{\partial}{\partial x^{1}}-\gamma \frac{\partial}{\partial x^{2}}+\frac{\partial}{\partial x^{3}} \right) ,
\end{displaymath}
where
\begin{displaymath}
\gamma :=\left( 1-\left( \frac{x^{3}}{\Arrowvert x \Arrowvert} \right) ^{\! 2} \right) x^{3}
\left\{ \left( 1-\left( \frac{x^{2}}{\Arrowvert x \Arrowvert} \right) ^{\! 2} \right) x^{2} \right\} ^{-1}.
\end{displaymath}
For $x^{2}=0$ we have $\Lambda (x)=0$ and
\begin{displaymath}
\ker d\phi _{x}=\spa \left( \frac{x^{1}}{x^{3}}\frac{\partial}{\partial x^{1}}+
\frac{\partial}{\partial x^{3}} \right) .
\end{displaymath}
In this case, we have that $\image d\phi _{x}\subseteq {\cal N}(TS^{2}_{1,1,1})$ but
$d\phi _{x}\neq 0$, i.e. we have case \emph{(b)} of Theorem \ref{bigguy}.

As we have
\begin{displaymath}
\frac{\partial ^{2}u}{(\partial x^{2})^{2}}=\frac{\partial ^{2}u}{(\partial x^{3})^{2}}=
\frac{2x^{2}x^{3}}{\Arrowvert x \Arrowvert ^{4}},
\end{displaymath}
then 
\begin{displaymath}
\overline{\tau}(\phi)=-\frac{\partial ^{2}u}{(\partial x^{2})^{2}}+
\frac{\partial ^{2}u}{(\partial x^{3})^{2}}=0
\end{displaymath}
so that $\phi$ is (generalized) harmonic.
By Theorem \ref{bigguy2}, the map $\phi$ is a (generalized) harmonic morphism.
\end{exam}
\section*{Aknowledgements}
I would like to thank Professor John C. Wood for his valuable support
and for commenting on drafts of this paper. I also thank the School
of Mathematics of the University of Leeds for the use of its facilities.
\bibliographystyle{amsplain}
\providecommand{\bysame}{\leavevmode\hbox to3em{\hrulefill}\thinspace}

\end{document}